%% file: paper.tex
\documentclass[final,leqno,onefignum,onetabnum]{siamltex1213}

\usepackage{amsmath}
\usepackage{amssymb}
\usepackage{bm}
\usepackage{booktabs}
\usepackage{color}
\usepackage{graphicx}
\usepackage{hyperref}
\usepackage{url}

\title{Interweaving PFASST and Parallel Multigrid}

\author{M.~L. Minion\footnotemark[1] 
\and R. Speck\footnotemark[2]~\footnotemark[3] 
\and M. Bolten\footnotemark[4]
\and M. Emmett\footnotemark[1]
\and D. Ruprecht\footnotemark[3]
}

\begin{document}

\maketitle
\slugger{sisc}{xxxx}{xx}{x}{x--x}

\renewcommand{\thefootnote}{\fnsymbol{footnote}}

\newcommand{\tvect}[1]{{\bf #1}}
\newtheorem{remark}{Remark}
\newcommand{\todo}[1]{\textcolor{red}{[TODO -- #1]}}
\newcommand{\note}[1]{\textcolor{blue}{[NOTE -- #1]}}

\footnotetext[1]{Center for Computational Sciences and Engineering, Lawrence Berkeley National Laboratory, USA}
\footnotetext[2]{J\"ulich Supercomputing Centre, Forschungszentrum J\"ulich, Germany}
\footnotetext[3]{Institute of Computational Science, Universit{\`a} della Svizzera italiana, Lugano, Switzerland}
\footnotetext[4]{Department of Mathematics, Bergische Universit\"at Wuppertal, Germany}

\renewcommand{\thefootnote}{\arabic{footnote}}
\def\1st {{$1^{\rm st}$~}}
\def\2nd {{$2^{\rm nd}$~}}
\def\4th {{$4^{\rm th}$~}}
\pagestyle{myheadings}
\thispagestyle{plain}
\markboth{M.~L. Minion, R. Speck, M. Bolten, M. Emmett, D. Ruprecht}{Interweaving PFASST and Parallel Multigrid} 

\begin{abstract}
  The parallel full approximation scheme in space and time (PFASST)
  introduced by Emmett and Minion in 2012 is an iterative strategy for
  the temporal parallelization of ODEs and discretized PDEs. As the
  name suggests, PFASST is similar in spirit to a space-time Full
  Approximation Scheme (FAS) 
  multigrid method performed over multiple time-steps in parallel. 
  However, since the original focus of
  PFASST has been on the performance of the method in terms of time
  parallelism, the solution of any spatial system arising from the use
  of implicit or semi-implicit temporal methods within PFASST have
  simply been assumed to be solved to some desired accuracy completely
  at each sub-step and each iteration by some unspecified procedure.
  It hence is natural to investigate how iterative solvers in the
  spatial dimensions can be interwoven with the PFASST iterations and
  whether this strategy leads to a more efficient overall approach.
  This paper presents an initial investigation on the relative
  performance of different strategies for coupling PFASST iterations
  with multigrid methods for the implicit treatment of diffusion terms
  in PDEs. In particular, we compare full accuracy multigrid solves at
  each sub-step with a small fixed number of multigrid V-cycles. This
  reduces the cost of each PFASST iteration at the possible expense of a
  corresponding increase in the number of PFASST iterations needed for
  convergence. Parallel efficiency of the resulting methods is
  explored through numerical examples.
\end{abstract}

\begin{keywords}
Parallel in time, PFASST, multigrid

\end{keywords}


\pagestyle{myheadings}
\thispagestyle{plain}

\input{defs}
\input{intro}

\input{sdc}

\input{mg_err}

\input{numerics}
\input{discussion}

\input{ack}

\bibliographystyle{siam}
\bibliography{sdc,Pint,Pint_Self,refs}

 \end{document}

%% file: defs.tex
\def\C{\mathbb C}
\def\R{\mathbb R}

\def\Dt {{\Delta t }}
\def\Dx {{\Delta x }}
\def\Nt {{N_t }}
\def\Nx {{N_x }}
\def\Dtm {{\Delta t_m }}
\def\dt {{\delta t }}
\def\tm {{t_m}}
\def\tmp {{t_{m+1}}}
\def\tn {{t_n}}
\def\tzero {{t_0}}
\def\tnp {{t_{n+1}}}
\def\np {{{n+1}}}

\def\Imat {{\bf I }}
\def\I {{ I }}
\def\Emat {{\bf E }}
\def\Pmat {{\bf P }}


\def\y {{y}}            
\def\yvec {{  \bf \y }}  
\def\Y {{Y}}            
\def\Yvec {{ \bf \Y }}   

\def\f   {{f}}          
\def\fvec   {{\bf{\f}}}  
\def\F   {{F}}          
\def\Fvec   {{\bf \F}}         

\def\yzero {{ \y_0 }}
\def\yveczero {{ \yvec_0 }}
\def\Yzero {{ \Y_0 }}
\def\Yveczero {{ \Yvec_0 }}

\def\yt {{{\y}(t)}}            
\def\yvect {{{\yvec}(t)}}  

\def\ytil {{\tilde{y}}}
\def\ytilm {{\tilde{y}_m}}
\def\ytilmp {{\tilde{y}_{m+1}}}
\def\ytilmk {{\tilde{y}^k_m}}
\def\ytilmpk {{\tilde{y}^k_{m+1}}}
\def\Ym {{\Y_m}}
\def\Ymp {{\Y_{m+1}}}
\def\Ymk {{\Y^k_m}}
\def\Ympk {{\Y^k_{m+1}}}
\def\Ymkp {{\Y^{k+1}_m}}
\def\Ympkp {{\Y^{k+1}_{m+1}}}
\def\Yt {{ \Y_p(t) }}

\def\Res {{R}}
\def\Restil {\tilde{R}}
\def\Ipmmp {{ \qmat^{m+1}_m }}
\def\Ipmmpvec {{ \Qmat^{m+1}_m }}

\def\mydel {{c}}
\def\mydelt {{\mydel(t)}}
\def\mydelkm {{\mydel^k_m}}

\def\myDelm {{\tilde{\mydel}_m}}
\def\myDelmp {{\tilde{\mydel}_{m+1}}}
\def\myDelkm {{\tilde{\mydel}^k_{m}}}
\def\myDelkpm {{\tilde{\mydel}^{k+1}_{m}}}
\def\myDelkpmp {{\tilde{\mydel}^{k+1}_{m+1}}}

\def\qmat {{ Q }}  
\def\qmatI {{ Q_I }} 
\def\Qmat {{\bf Q }} 
\def\QmatI {{\bf Q_I }} 

\def\ellp {{\ell+1}}
\def\QmattriC {{\bf Q}^\ellp_{\bf \triangle}}
\def\Qmattril {{\bf Q}^\ell_{\bf \triangle}}
\def\xvec {{\bf x}}

\def\Smat {{\bf S }}
\def\SmatE {{\bf S_E }}
\def\SmatI {{\bf S_I }}
\def\qvec {{\bf q}}
\def\tvec {{ \bf t }}
\def\zerovec {{ \bf 0 }}

\def\Amat {{\bf A }}
\def\cvec {{\bf c}}
\def\Cvec {{\bf C}}

\def\bigQ {{\bf \Qmat }}
\def\bigQtri {{\bf \bar{Q}_{\triangle}}}
\def\bigQtrilp {{\bf \bar{Q}}^\ellp_{\triangle}}
\def\bigQI {{\bf \bar{Q}_{I}}}
\def\bigQE {{\bf \bar{Q}_{E}}}
\def\bigQtriC {{\bf \bar{Q}}^C_{\bf \triangle}}

\def\bigA {{\bf A }}
\def\A {{A }}
\def\M {{\bf M }}
\def\P {{\bf P }}
\def\C {{\bf C }}
\def\S {{\bf S }}
\def\K {{\bf K }}
\def\D {{\bf D }}
\def\L {{\bf L }}
\def\U {{\bf U }}
\def\bigAE {{\bf \bar{A} }_E}
\def\bigAI {{\bf \bar{A} }_I}
\def\bigI {{\bf I }}

\def\bigM {{\bf M }}
\def\bigK {{\bf K }}


\def\Fvec {{ \bf F }}
\def\FvecE {{ \bf F }_E}
\def\FvecI {{ \bf F }_I}

\def\Onevec {{\bf 1 }}
\def\onevec {{\bf 1 }}
\def\Smattri {{\bf S_{\triangle} }}
\def\Qmattri {{\bf Q_{\triangle} }}
\def\PmatE {{\bf P_E }}
\def\PmatI {{\bf P_I }}
\def\Pmattri {{\bf P_\triangle }}
\def\Pmattril {{\bf P}^\ell_{\bf \triangle }}
\def\PmattriC {{\bf P}^\ellp_{\bf \triangle }}


\def\Lvec {{ \bf L }}
\def\Ibar {{\bf \bar{I} }}
\def\Pbar {{\bf \bar{P} }}
\def\Ybar  {{\bf \bar{Y} }}
\def\Qmatbar {{\bf \bar{Q} }}
\def\Qmatbartri {{\bf \bar{Q}_{\triangle} }}
\def\QmatbarE {{\bf \bar{Q}_E }}
\def\QmatbarI {{\bf \bar{Q}_I }}

\def\Yvecn {{ \Yvec_n }}
\def\Ymn {{ \Y_{m,n} }}
\def\lamdt {{\lambda \Dt}}
\def\lamDt {{\lamdt}}
\def\qzvec {{\bf q_0}}
\def\Ybarbar {{\bf \bar{\bar{Y}}}}

\def\ellp {{\ell +1 }}
\def\ellm {{\ell -1 }}
\def\bvec {{\bf{b}}}
\def\ResT {{\bf R}_\ell^\ellp}
\def\ResTaulp {{\bf R}^\ellp_\tau}
\def\ResTaul {{\bf R}^\ell_\tau}
\def\ProT {{\bf T}_\ellp^\ell}
\def\bigResT {{\bf \bar{R}_\ell^\ellp}}
\def\bigResTaulp {{\bf \bar{R}^\ellp_\tau}}
\def\bigProT {{\bf \bar{T}_\ellp^\ell}}
\def\Sweepl {{\bf \bar{S}^\ell}}
\def\Solvelp {{\bf \bar{S}^{-\ellp}}}
\def\Residl {{\bf \bar{R}^\ell}}

%% file: intro.tex
\section{Introduction}

The past decade has seen a growing interest in the development of
parallel methods for temporal integration of ordinary differential
equations (ODEs), particularly in the context of temporal strategies
for partial differential equations (PDEs).  One factor fueling this
interest is related to the evolution of supercomputers during this
time.  Since the end of the exponential increase in individual
processors speeds, increases in supercomputer speeds have been mostly
due to increases in the number of computational cores, and current
projections suggest that the first exaflop computer will contain on
the order of a billion cores~\cite{Geist2010}.  The implication of this
trend is that increasing concurrency in algorithms is essential, and
in the case of time-dependent PDE simulations, the use of space-time
parallelism is an attractive option.

Time-parallel methods have a long history dating back at least to the
work of Nievergelt~\cite{Nievergelt1964}.  In the context of
space-time multigrid, Hackbusch noted already in 1984 that relaxation
operators in parabolic multigrid can be employed on multiple time
steps simultaneously~\cite{Hackbusch1984}.  The 1997 review article by
Burrage~\cite{Burrage1997} provides a summary of early work on the
subject.  More recently, the parareal method proposed in
2001~\cite{LionsEtAl2001} has renewed interest in temporal
parallelization methods.  In 2012 the parallel full approximation
scheme in space and time (PFASST) was introduced by Emmett and
Minion~\cite{EmmettMinion2012}, and performance results for PFASST 
using space-time parallelization with hundreds of thousands of cores
can be found in~\cite{SpeckEtAl2012,RuprechtEtAl2013_SC}.  

The PFASST algorithm is based on a type of deferred corrections strategy for
ODEs~\cite{DuttEtAl2000}, with corrections being applied on multiple
time steps in parallel.  As such, there are similarities between
parareal and PFASST (see~\cite{MinionEtAl2008,Minion2010}).
On the other hand, the parallel efficiency of PFASST depends on the
construction of a hierarchy of space-time discretizations, hence there
are also similarities between PFASST and space-time multigrid methods.
However, in the original papers on PFASST, the solution of any spatial systems
due to implicit time-stepping was assumed to be found to full
precision since the interest was on the temporal accuracy and
efficiency of the methods.  From this point of view, PFASST is an
iterative solver in the time direction but not in the spatial
dimensions.  This is, in a sense, orthogonal to the traditional use of
multigrid solvers within PDE methods, where multigrid is used to
iteratively solve spatial equations and the time direction is not
iterative.  One of the main goals of this paper is to investigate the
use of both spatial and temporal iterative methods utilizing PFASST
and multigrid. 

To be more specific, the iterative strategy within the PFASST algorithm is
derived from the method of Spectral Deferred Corrections~\cite{DuttEtAl2000},
which is a variant of the defect and deferred correction
methods developed in the 1960s~\cite{Zadunaisky1966a,Pereyra1966,stetter:1974,bohmerHemkerStetter:1984}.  One advantage of SDC methods
is that it is straightforward to construct methods of very high order of
accuracy by iteratively applying low-order methods to a series of correction
equations.  This flexibility has been exploited to construct higher-order 
semi-implicit (IMEX) and multi-implicit SDC  methods~\cite{Minion2003,BourliouxEtAl2003,LaytonMinion2004}, as well
as multirate SDC methods~\cite{BouzarthMinion2010}.  Such methods are very difficult
to construct using traditional Runge-Kutta or linear-multistep approaches.

The main disadvantage of SDC methods is that they have a high cost per time
step in the sense of the number of function evaluations (explicit or implicit)
required per time step.  When high accuracy is desired, this cost is generally
offset by the use of relatively large time steps compared to lower-order 
methods for the same given accuracy.  Nevertheless, approaches to 
reducing the cost of each SDC iteration have been proposed, including those
generally referred to as {\it multi-level} SDC methods (MLSDC).  The
main idea in MLSDC is to perform some SDC iterations on coarser discretizations
of the problem, and methods that coarsen the temporal discretization, the spatial
discretization, and the order of the spatial approximation have recently been
investigated~\cite{SpeckEtAl2014_BIT}.  SDC iterations are then performed on the hierarchy
of levels in much the same way as V-cycles in traditional multigrid.  The PFASST
algorithm can be considered  a time parallel version of MLSDC.

Using implicit SDC methods (as well as popular methods like backward Euler, trapezoid rule,
Diagonally Implicit Runge-Kutta, or Backward Difference Formula)
for grid-based PDEs simulations requires the solution
of implicit equations in the spatial domain. If an iterative solver is used for 
these equations in the context of SDC, one advantageous result is that the initial
guess for each solve becomes better as the SDC iterations converge. This raises
the possibility of reducing the cost of SDC or MLSDC iterations further by limiting
the number of spatial iterations used for each implicit solve rather than requiring
each to be done to some specified tolerance.  
A recent paper explores this possibility in the context of SDC methods~\cite{SpeckEtAl2014_DDM2013}. 
We refer to variants in which spatial solves are done
only incompletely with a prepended capital ``I'' (ISDC, IMLSDC, IPFASST), resulting in the
interweaving of spatial and temporal iterations. 
This paper adopts the ISDC strategy for MLSDC and PFASST and explores the performance of both IMLSDC and IPFASST.

As a starting point, we focus here on the the linear test problem
resulting from a finite difference approximation of the heat equation.
Variants of this example are considered in early papers on space-time
multigrid~\cite{HortonKnirsch1992,HortonVandewalle1995}, waveform
relaxation~\cite{VandewalleHorton1995,Gander1998}, block-parallel
methods~\cite{AmodioBrugnano2009}, as well as more recent papers on
parallel space-time multigrid~\cite{FriedhoffEtAl2013,FalgoutEtAl2014_MGRIT}.  In all the
cases above, a second-order finite-difference approximation is
employed in space and first- or second-order methods used in time.
Although this is not the optimal setting for the PFASST algorithm, we
present results using second-order space-time discretization to
compare with other published results.  We also show how, for a given
accuracy, using fourth-order methods in space and/or time results in
significant computational 
savings compared to second-order methods.

The remainder of this paper is organized as follows.  We first present
the SDC algorithm for linear problems in a compact notation that can
be interpreted as matrix iterations and highlight the similarities
between SDC iterations and classical relaxation schemes. We then
discuss how IMLSDC and IPFASST are constructed. 
In Section \ref{sect:linear}, we provide the relevant analysis
of SDC methods as a relaxation operator for parabolic problems.
In Section \ref{sect:numerics} we examine the
scaling of the IPFASST method in terms of problem size and number of
parallel time steps, consider the effect of limiting the number of
multigrid V-cycles per implicit solve, and provide numerical examples comparing
PFASST with IPFASST.  Finally, in Section \ref{sect:strong}
we provide timing comparisons for
three-dimensional examples scaling to hundreds of thousands of cores.

%% file: sdc.tex
\section{Collocation and Spectral Deferred Corrections}

This section briefly describes  the methods used later.
In Section~\ref{subsect:picard}, the collocation formulation for initial value
problems is first reviewed.  SDC as an iterative solver for a collocation
solution is discussed briefly in Section~\ref{subsect:sdc}.  A compact
notation of SDC for linear problems is given, and its interpretation as
relaxation is discussed.  The extension of SDC to multi-level SDC
(MLSDC) is outlined in  Section~\ref{subsec:mlsdc}, including a number of
strategies to coarsen the representation of the problem 
in order to reduce 
the overall cost.  
The possibility of 
solving linear systems approximately on all
levels leads to so-called "inexact spectral deferred corrections"
(ISDC) and their multi-level counterparts IMLSDC and IPFASST.

\subsection{The Collocation Formulation} \label{subsect:picard}

This paper concerns methods for 
the solution of initial value ODEs, particularly those
arising from the spatial discretization of PDEs through
the method of lines technique.
To set notation, consider the scalar ODE of the form
\begin{eqnarray}
  \y'(t) &  = & \f(t,\y(t))  \label{eq:ODE} \\
  \y(t_0) & = & \y_0,
\end{eqnarray}
where $y(t), \y_0 \in \C$ and $f: \R \times \C \rightarrow \C$.
Similarly for systems of equations
\begin{eqnarray}
  \yvec'(t) &  = & \fvec(t,\yvect)  \label{eq:sysODE} \\
  \yvec(t_0) & = & \yvec_0,
\end{eqnarray}
where $\yvec(t), \yvec_0 \in \C^D$ and $\fvec: \R \times \C^D
\rightarrow \C^D$.   An equivalent form is given by the Picard
equation, which for \eqref{eq:sysODE} is
\begin{equation}
  \yvect =  \yvec_0 + \int^t_{t_0} \fvec(\tau,\yvec(\tau)) d\tau. \label{eq:Picard}
\end{equation}

Since the goal here is to describe numerical
methods for a given time step, consider the time interval $t \in
[\tn,\tnp]$ with $\Dt = \tnp-\tn$.  
Define the set of points $t_m$ for   $m=0, \hdots, M$ to
be quadrature nodes scaled to $[\tn,\tnp]$, 
so that $\tn=t_0 < t_1 < t_2 < \hdots < t_M = \tnp$.  
Given a scalar function $g(t)$, define the 
approximations to the integral of $g(t)$ over the intervals 
$[\tzero,t_m]$ by choosing the coefficients $q_{m,i}$ such that
\begin{displaymath}
\int_{t_0}^\tm g(\tau) d\tau \approx \Dt \sum_{i = 0}^M ~q_{m,i} g(t_i).
\end{displaymath}
The coefficients $q_{m,i}$ that give the highest order
of accuracy given the points $t_m$ 
can be derived using standard techniques.  The quadrature rules
used here have the property that $q_{m,0} \equiv 0$ (see 
\cite{LaytonMinion2005} for a discussion of different choices
of quadrature rules for SDC methods).
Let the matrix $\qmat$ of size $(M+1) \times
(M+1)$ 
be composed of the coefficients $q_{m,i}$.  The first row of
$\qmat$  contains only zeros.  $\qmat$ will be referred to here as the integration matrix.

Consider for the time being scalar equations, the integration matrix $\qmat$  
can be used to discretize \eqref{eq:Picard} directly
over the interval $[\tn,\tnp]$. Let
$y_0=y_n \approx y(t_n)$ 
be the initial
condition, and define $\y_m \approx \y(t_m)$ by 
\begin{equation} 
	\label{eq:collocation}
	\y_m = \y_n + \Dt \sum_{i = 0}^M  q_{m,i} \f(t_i,\y_i)  \hbox{\hspace{0.5cm} for } m = 0, \ldots, M.
\end{equation}
We can write this equation in a compact form by denoting the vector of unknowns by $\Y = [\y_0 \ldots \y_M]^T$, and the vector of function values 
$\F(\Y) = [\f(\y_0,t_0) \ldots \f(\y_M,t_M)]^T$.
Furthermore, let  $\Yzero$ be the $(M+1) \times 1$ column vector with each entry equal to $\yzero$.
Then \eqref{eq:collocation} is equivalent to 
 \begin{equation} 
 	\label{eq:collocationM}
 	\Y    = \Yzero +  \Dt \qmat \F(\Y).
\end{equation}
This coupled (generally nonlinear) equation can 
be solved directly for the values $\Y$,
resulting in a collocation scheme for the ODE. 

It should also be noted 
that
the collocation form 
\eqref{eq:collocationM} corresponds to a fully implicit Runge-Kutta  (IRK)
method given by the Butcher tableau
\begin{equation}  \label{eq:BT}
\begin{tabular}{ c | c  }
$\cvec$ & $\Amat$   \\
\hline
        & $\bvec$ \\
\end{tabular},
\end{equation}
where $\cvec$ are the nodes $t_m$ scaled to the unit interval,
$\Amat=\qmat$, 
and the vector $\bvec$ corresponds 
to the last row of $\qmat$.

For the simplest linear scalar equation where $\f(y,t) = \lambda y$ in \eqref{eq:ODE}, the collocation formulation given in \eqref{eq:collocationM} becomes 
 \begin{equation} \label{eq:collocationML}
 \Y    = \Yzero +  \lamdt \qmat  \Y,
\end{equation}
or in more standard matrix form
 \begin{equation} \label{eq:collocationMLII}
(\I_{M+1} - \lamdt \qmat) \Y    = \Yzero,
\end{equation}
with $\I_{M+1}$ being the $(M+1) \times (M+1)$ identity matrix.
In the case of a linear system of equations 
$\fvec(t,\yvec) = \A \yvec$, where $\A$ is a $D \times D$ matrix,
the Picard equation becomes
\begin{equation}
  \yvec(t) =  \yvec_n + \int^t_{t_n} \A \yvec(\tau) d\tau. \label{eq:LinPicard}
\end{equation}
where the integration is done term by term.
To discretize \eqref{eq:LinPicard}, the solution is the vector 
\begin{equation}
\Yvec = [\y^1_0 \ldots \y^D_0 \ldots \y^1_M \ldots \y^D_M]^T
\end{equation}
of length $D(M+1)$.
Here the subscripts $m$ correspond to the quadrature nodes in time 
and superscripts $j$ correspond to the component of the solution vector.

To apply the numerical quadrature matrices to $\Yvec$,
let $\Qmat = \qmat \otimes \I_D$ where $\I_D$ is the $D \times D$ identity
matrix.   Likewise let $\Amat = \I_{M+1} \otimes \A$, and $\Yveczero$ be
the vector of length $(M+1)D$ consisting of $M$ copies of the $\yveczero$. 
Then the analogous form of \eqref{eq:collocationM}  is
\begin{equation} \label{eq:bigPmatk_1}
\Yvec = \Yveczero + \Dt \bigQ \bigA \Yvec
\end{equation}
or
\begin{equation} \label{eq:bigPmatk_2}
(\Imat -  \Dt \bigQ \bigA) \Yvec = \Yveczero ,
\end{equation}
where $\Imat = \I_{M+1} \otimes \I_D$.


\subsection{Single level SDC}\label{subsect:sdc}
SDC can be understood as a preconditioned fixed point iteration to
solve~\eqref{eq:bigPmatk_2}, which avoids treating the full system
directly by computing a series of corrections  node by node (see
e.g.~\cite{HuangEtAl2006,Weiser2013}).  Note that originally SDC was
introduced in a different way, as a variant of earlier deferred and 
defect correction schemes \cite{Zadunaisky1966a,Pereyra1966}
designed to achieve a fixed order of accuracy for a fixed number
of correction sweeps~\cite{DuttEtAl2000}.  The multi-level SDC
methods (and PFASST) described below move away from the idea
of a fixed number of iterations in favor of the convergence toward
the collocation (or IRK) solution.

Here we provide a short review of the SDC method. For more details see
\cite{DuttEtAl2000,Minion2003}.  Let the superscript $k$ denote the
numerical approximation at the $k$th SDC iteration.  Using
backward-Euler as the base method, one iteration of SDC can be written
as a sweep through the quadrature nodes, successively updating the
solutions.  In the case of scalar equations, this takes the form
\begin{equation}\label{sdc-be}
\y^{k+1}_{m+1} = \y^{k+1}_m  +  \Dtm  \left(\f(\tmp,\y^{k+1}_{m+1}) 
                                          -   \f(\tmp,\y^{k}_{m+1})\right) +  \Ipmmp(\Y^{k})
\end{equation}
for $m=0, \ldots, M-1$ where
\begin{equation}
 \Ipmmp(\Y^{k})   \approx  \int^\tmp_\tm  \f(\tau,\y^k(\tau)) d\tau.
\end{equation}
The values of  $\Ipmmp(\Y^{k})$ can easily be constructed using the integration matrix
$\Qmat$
\begin{equation}
 \Ipmmp(\Y^{k})  =  \Dt \sum_{i = 0}^M ~(q_{m+1,i}-q_{m,i}) \f(t_i,\y^k_i).
\end{equation}
Here we consider only implicit SDC methods, even though one particularly attractive feature of SDC is the flexibility to use different base methods in order to create e.g. high-order implicit-explicit or multi-rate methods~\cite{Minion2003,BourliouxEtAl2003,HagstromZhou2006,BouzarthMinion2010}.

\subsubsection{Compact notation} 
We will now present a compact notation of the SDC iterations.
Given the points $t_m \in [\tn,\tnp]$ as discussed above,
let $\Dt_m = t_m-t_{m-1}$, and let $\gamma_m = \Dt_m/\Dt$,
where again $\Dt = t_{n+1}-t_{n}$.
We begin by defining the lower-order integration matrix
that has the same dimensions as $\qmat$
 \begin{equation}\label{eq:Qbe}    
 \qmatI  = \left[ \begin{array}{ccccc}
  0 & 0  & \cdots & 0 & 0 \\
  0 & \gamma_1  & \cdots & 0 & 0 \\
  \cdot & \cdot & \cdots & 0 & 0 \\
  0 & \gamma_1  & \cdots & \gamma_{M-1} & 0 \\
  0 & \gamma_1  & \cdots & \gamma_{M-1} & \gamma_{M}
 \end{array} \right].
 \end{equation}
Using the same notation as in Section \ref{subsect:picard}, 
$M$ steps of backward Euler for \eqref{eq:ODE} 
can be written as
\begin{equation} \label{eq:bE}
 \Y  = \Yzero + \Dt \qmatI \F(\Y).
\end{equation}
This differs from the collocation formulation \eqref{eq:collocationM}
only in that $\qmatI$ replaces $\qmat$.

Since $\qmatI$ is lower triangular with non-zero diagonal entries, \eqref{eq:bE}
can be solved sequentially for the values $y_m$, where each value
requires an implicit equation to be solved of the form
\begin{displaymath}
 \y_{m+1}  = \y_m + \Dt_m  \f(t_{m+1}, \y_{m+1}).
\end{displaymath}
Using the matrix-vector  notation above, 
one SDC sweep for \eqref{eq:sysODE} can be compactly written as
 \begin{equation} \label{eq:SDCcompact}
 \Yvec^{k+1}    = \Yveczero
+ \Dt \QmatI ( \Fvec(\Yvec^{k+1})-\Fvec(\Yvec^k))  + \Dt \Qmat \Fvec(\Yvec^k)
\end{equation}
or 
 \begin{equation} \label{eq:SDCcompactII}
 \Yvec^{k+1}    = \Yveczero
+ \Dt  \QmatI \Fvec(\Yvec^{k+1})  + \Dt (\Qmat-\QmatI) \Fvec(\Yvec^k).
\end{equation}

\subsubsection{SDC for Linear Problems} \label{subsub:compact}

The compact formulas derived above can be recast as 
matrix-vector operations  when the governing equation 
is linear. In the case of the linear system, $\yvec'= 
\A \yvec$,  \eqref{eq:SDCcompact}
simplifies to 
 \begin{equation} \label{eq:SDCcompactL}
 \Yvec^{k+1}    = \Yveczero
+ \Dt \QmatI ( \Amat \Yvec^{k+1}- \Amat \Yvec^k )  + \Qmat \Amat \Yvec^k
\end{equation}
or 
 \begin{equation} \label{eq:SDCcompactLII}
\left( \Imat -  \QmatI \Amat \right) \Yvec^{k+1}    = \Yveczero
+  \Dt \left(\Qmat-\QmatI\right) \Amat \Yvec^k.
\end{equation}
As in the case of backward Euler, this system of equations can be
solved by sub-stepping, requiring the solution of 
 \begin{equation} \label{eq:SDCcompactLIII}
( \I - \Dt_m \A ) \yvec^{k+1}_{m+1}    = \yveczero-  \Dt_m \A\yvec^{k}_{m+1}
+ \Ipmmpvec(\Yvec^{k}).
\end{equation}
This equation involves the inversion of the same operator that arises
from a backward Euler method with a modified right-hand side.  Note 
that as the SDC iterations converge in 
$k$, an increasing good approximation to this solution is provided
by $\yvec^{k}_{m+1}$.  This fact has two relevant implications when
iterative methods are employed to solve the system.  First, when
considering the total cost of one time-step of SDC, the reduced
cost of the implicit solves as $k$ increases should be taken into
consideration.  Second, and more relevant for this paper, is that
an iterative method need not solve \eqref{eq:SDCcompactLIII} to
full precision at each iteration.  Instead, a fixed number of iterations
could be performed during each SDC iteration, or each could be
done so that the residual is reduced by some set tolerance.
In the numerical studies presented in Section \ref{sect:numerics},
multigrid methods are employed for solving \eqref{eq:SDCcompactLIII}
and we investigate limiting the number of multigrid V-cycles 
instead of solving to full precision.

Finally, note that in the SDC iterations, it is straightforward
to compute the residual in terms of the solution to 
\eqref{eq:bigPmatk_1}, specifically
\begin{equation}
	\label{eq:resid}
	\tvect{r}^{k} := \Yveczero - (\bigI - \Dt \bigQ \bigA)\Yvec^{k}.
\end{equation}
In Section~\ref{sect:numerics}, the norm of $\tvect{r}^k$ 
is used to monitor convergence of the different choices of methods.

\subsubsection{SDC as a relaxation}\label{subsubsect:sdc_relax}
For a linear problem, it is possible to write an SDC sweep described 
above as a relaxation  operator applied to the linear collocation
equation 
\begin{equation} \label{eq:bigPicard_1}
(\bigI - \Dt \bigQ \bigA)\Yvec = \Yveczero.
\end{equation}
As in classical iteration methods for linear systems, we look for a decomposition
of the matrix $(\bigI - \Dt \bigQ \bigA) = \bigM-\bigK$ such that $\bigM$ is 
relatively less expensive to invert than $(\bigI - \Dt \bigQ \bigA)$.  Then
a classical relaxation scheme based on this splitting would be
\begin{equation} \label{eq:bigrelax}
\Yvec^{k+1} =   \bigM^{-1}\bigK \Yvec^k + \bigM^{-1} \Yveczero.
\end{equation}
Choosing
\begin{equation} \label{eq:bigPicard_2}
\bigM = \bigI - \Dt \QmatI \bigA \hbox{ \hspace{0.25cm} and  \hspace{0.25cm}}  \bigK = \Dt (\Qmat-\QmatI) \bigA 
\end{equation}
produces the SDC sweep as given by \eqref{eq:SDCcompactLII}, hence
inverting $\bigM$ can be done by sub-stepping on the SDC nodes solving
the appropriate version of \eqref{eq:SDCcompactLIII}.

In the next section we carry this analogy further by introducing
variants of SDC that utilize multiple levels of resolution as in
classical multigrid methods.  

\subsection{Multi-level SDC Methods (MLSDC)}\label{subsec:mlsdc}
\begin{figure}[t]
	\centering
	\includegraphics[height=0.15\textheight]{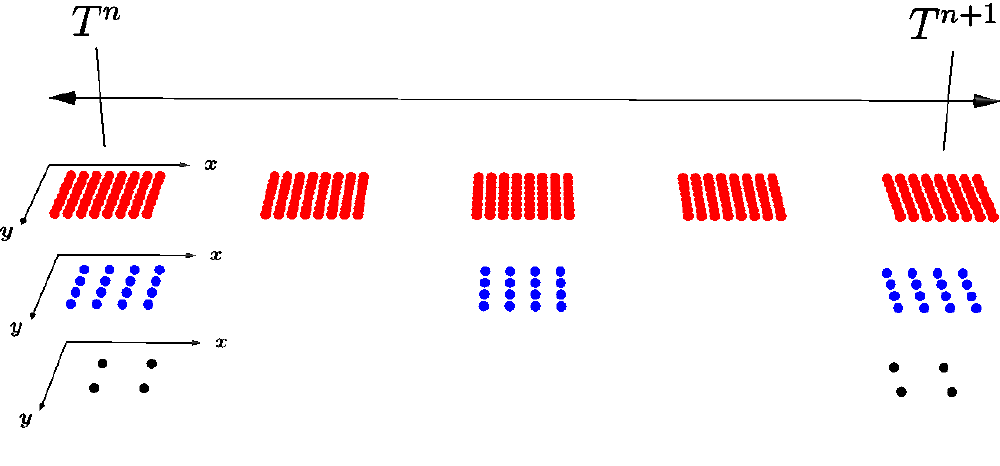}
	\caption{A MLSDC level hierarchy over a time-step $[T^n, T^{n+1}]$. 
The upper fine level (red) uses five quadrature nodes, the medium level (blue) three and the bottom coarsest level (black) two. Also sketched is the possibility of using a spatial discretization with fewer degrees of freedom on the coarser levels: Depicted is an $8 \times 8$ point mesh on the fine level, a $4 \times 4$ point mesh in the middle, and a $2 \times 2$ point mesh on the coarse level. A detailed description of MLSDC can be found in~\cite{SpeckEtAl2014_BIT}. Figure reprinted from~\cite{RuprechtEtAl2013_SC}.}\label{fig:mlsdc}
\end{figure}

The goal of MLSDC methods is that by introducing a hierarchy of levels 
from fine to coarse, some of the expensive fine level 
correction 
sweeps can be replaced with sweeps on coarser levels, so that
the runtime required for convergence of SDC iterations is reduced (see recent results in~\cite{SpeckEtAl2014_BIT}).  
An space-time 
FAS term is employed when forming coarser level 
equations, which is 
the difference between the coarsening of the
latest fine function values and the coarse level function
applied to a coarsening of the latest fine solution.
The result of including the FAS term is that the 
accuracy on the coarse level approaches that of the fine
level as the MLSDC iterations converge (see for example \cite{multigridtutorial}).
Different strategies for
reducing  the cost of SDC sweeps on the coarser levels are explored
in \cite{SpeckEtAl2014_BIT}, including using a lower-order spatial
discretization, using a spatial mesh with fewer points, 
or performing only incomplete implicit solves.  
The structure of an MLSDC level hierarchy is sketched in
Figure~\ref{fig:mlsdc}. For a detailed explanation of the method
including performance results we refer to~\cite{SpeckEtAl2014_BIT}.
MLSDC is the basic building block for the time-parallel PFASST method
summarized in Section~\ref{subsec:pfasst}.

\subsubsection{Inexact MLSDC (IMLSDC)}
In~\cite{SpeckEtAl2014_DDM2013}, another approach to reducing
the cost of SDC methods is discussed, 
so called "inexact spectral deferred corrections" or ISDC.
Here, implicit solves in SDC sweeps~\eqref{sdc-be} are computed only
approximately, e.g. by a fixed small number of V-cycles of 
multigrid.  While this strategy makes each SDC sweep less expensive, it
can increase the number of SDC iterations required to 
reduce the residual to a given threshold.  If the increase in
iterations is not too large, ISDC can in total be faster than SDC.
In~\cite{SpeckEtAl2014_DDM2013} it is demonstrated that ISDC can save
up to $50 \%$ of the V-cycles required by regular SDC while converging
to the same tolerance.  It is also shown that the reason this strategy
works is the increasingly accurate initial guesses provided for the
iterative solver by SDC (see the discussion toward the end of 
Section~\ref{subsub:compact}).

In regular MLSDC, doing incomplete solves on the coarse levels is
already discussed and explored in~\cite{SpeckEtAl2014_BIT} as a
means to reduce the cost of coarse level sweeps.  However, the ISDC
strategy can also be incorporated into MLSDC, resulting in an IMLSDC
method. Here, all implicit solves are done only approximately, including
the finest level.  Ideally, IMLSDC would save on runtime compared to
ISDC, just as MLSDC does compared to SDC.  However, the results
presented in Section~\ref{sect:numerics} suggest that this is not
necessarily the case.  However, MLSDC serves as the basis of the 
parallel-in-time 
method 
PFASST, and hence IMLSDC could be used 
in parallel as well. In Section~\ref{sect:numerics}
we demonstrate that IPFASST  (see Section~\ref{subsubsect:ipfasst}) can provide
a significant reduction of runtime by exploiting temporal concurrency.
A detailed description of MLSDC including pseudo code can be 
found in~\cite{SpeckEtAl2014_BIT}, and we refer the reader there for details.
%
%


\subsection{The Parallel Full Approximation Scheme in Space and Time}\label{subsec:pfasst}
The parallel full approximation scheme in space and time (PFASST),
introduced in~\cite{Minion2010,EmmettMinion2012} is an iterative time-parallel method for PDEs that has similarities to both
the parareal method and space time multigrid methods.  In PFASST each
time step is assigned its own processor or, if combined with spatial
parallelization, its own communicator (for a more detailed explanation
of the latter case, see~\cite{SpeckEtAl2012}).  PFASST can be
considered a time-parallel extension of MLSDC, where after an
initialization procedure, MLSDC iterations are performed on multiple
time steps in parallel with updates to initial conditions being passed
between processors as each SDC sweep is completed. 

\begin{figure}[t]
	\centering
	\includegraphics[height=0.15\textheight]{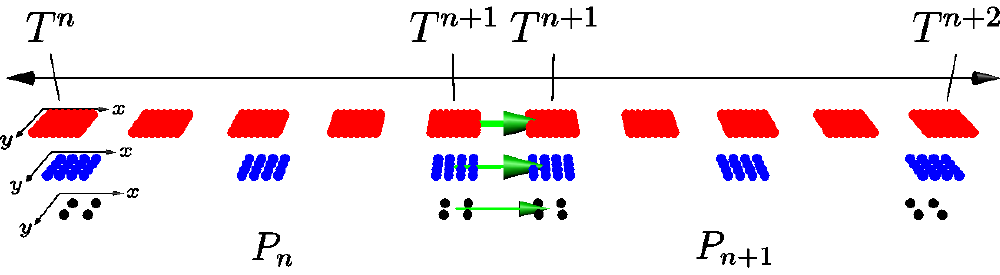}
	\caption{PFASST performs MLSDC iterations concurrently on
          multiple time steps. For simplicity, only two time steps
          $[T^n, T^{n+1}]$ and $[T^{n+1}, T^{n+2}]$ are sketched
          here. After each sweep, an updated initial value on the
          corresponding level is 
          sent forward to the processor
          handling the next time step. In the setup above, processor
          $P_n$ handles time step $[T^n, T^{n+1}]$ and sends forward
          updates to processor $P_{n+1}$ which handles $[T^{n+1},
          T^{n+2}]$. Figure reprinted
          from~\cite{RuprechtEtAl2013_SC}.}\label{fig:pfasst}
\end{figure}

PFASST typically
starts by distributing the initial value over all time ranks as an
initial guess which is then refined by a number of sweeps on the
coarse level, where processors handling time steps later in time do
more sweeps (this is usually referred to as predictor phase).
After completing the predictor phase, each processor
starts with its own MLSDC iterations while, after each sweep, sending
forward an updated initial value for the current level to the
processor handling the next time step, cf. Figure~\ref{fig:pfasst}.
Blocking communication is required on the coarsest level only, so that
minimal synchronicity between the different MLSDC iterations is
required, see~\cite{EmmettMinion2014_DDM}.  Benchmarks for PFASST in
large-scale parallel simulations illustrate how PFASST can extend
strong scaling limits~\cite{SpeckEtAl2012} or improve utilization of
large parallel installations in comparison to codes utilizing only
spatial parallelism~\cite{RuprechtEtAl2013_SC}.

A detailed description of PFASST including a sketch of the algorithm in pseudo code can be found in~\cite{EmmettMinion2012}, and we refer the reader there for more details.
\subsubsection{IPFASST}\label{subsubsect:ipfasst}
Just as incomplete implicit solves  can be used in
MLSDC, they can also be used in PFASST.  Essentially, each processor
now performs IMLSDC iterations instead of iterations with full solves
on the fine level.  Apart from that, IPFASST 
proceeds the same
as PFASST. In Section~\ref{sect:numerics}, performance of IPFASST 
will be studied through numerical examples.

%% file: mg_err.tex
\section{Analysis of SDC on the Linear Model Problem} \label{sect:linear}

In Section \ref{subsubsect:sdc_relax}, the analogy between
a single SDC sweep and a classical relaxation scheme
derived from a splitting of the linear system is presented.
In this section we examine the effect of a single
SDC sweep on the scalar linear model problem  
\begin{eqnarray}
y' & = & \lambda y \\
y(0) & = & 1. 
\end{eqnarray}
In the context of parabolic problems, the relative set
of $\lambda$ is the negative real axis, and we are
interested in how one sweep of SDC reduces the 
error in the discrete approximation.

In this case,  the compact form of the SDC sweep from
\eqref{eq:SDCcompactLII} becomes
 \begin{equation} \label{eq:lamdtI}
(\I - \lamDt \qmatI) \Y^{k+1}    =\Yzero 
+ \lamdt ( \qmat-\qmatI  )\Y^k.
\end{equation}
The corresponding solution to the collocation equation satisfies
 \begin{equation} \label{eq:lamdtII}
(\I - \lamDt \qmat) \Y    = \Yzero,
\end{equation}
which implies
 \begin{equation}   \label{eq:lamdtIII}
(\I - \lamDt \qmatI) \Y    = \Yzero
+ \lamdt ( \qmat-\qmatI)\Y.
\end{equation}
In words,  the SDC sweep acts as the identity operator for
the exact solution $\Y$.

Subtracting equation \eqref{eq:lamdtIII} from \eqref{eq:lamdtI}
gives
 \begin{equation} \label{eq:delsweep}
(\I - \lamDt \qmatI )(\Y^{k+1}-\Y)  = 
\lamDt (\qmat-\qmatI )(\Y^k-\Y)
\end{equation}
or
 \begin{equation} \label{eq:delsweepII}
\Y^{k+1}-\Y  = (\I - \lamDt \qmatI )^{-1}\lamDt(\qmat-\qmatI )(\Y^k-\Y).
\end{equation}
Hence we can study the converge properties of SDC methods for the 
scalar model problem by examining the largest eigenvalue of the matrix
\begin{equation} \label{eq:SDCMAT}
   (\I - \lamDt \qmatI )^{-1}\lamDt(\qmat-\qmatI ).
\end{equation}
This matrix depends in general on the form of the sub-stepping
encapsulated in the approximate quadrature matrix (here
backward Euler in 
$\qmatI$), 
the number and type of quadrature nodes,
and the product $\lamdt$.  Here we examine the cases corresponding
to the second- and fourth-order methods used in the numerical
results in Section \ref{sect:numerics}.  These correspond
to uniform quadrature nodes and a quadrature rule which does
not use the left-hand endpoint in the quadrature rule
(see \cite{LaytonMinion2005} for a discussion of different
quadrature rules). These methods are used together in the 
fourth-order PFASST example as well, where the second-order
method is the time coarsened version of the fourth-order method.
This choice requires 2 implicit sub-steps
for second order and 4 for fourth order and hence is not
``spectral'' in the sense of using Gaussian quadrature rules.
This choice does however provide good damping factors for
low-order methods\footnote{A recent paper by Weiser \cite{Weiser2013}
studies the optimization of more general sub-stepping rules for SDC methods.
These ideas are not pursued used here despite their apparent promise.}.
\begin{figure}[t] \label{fig:damping}
\centering
 \includegraphics[width=0.7\textwidth]{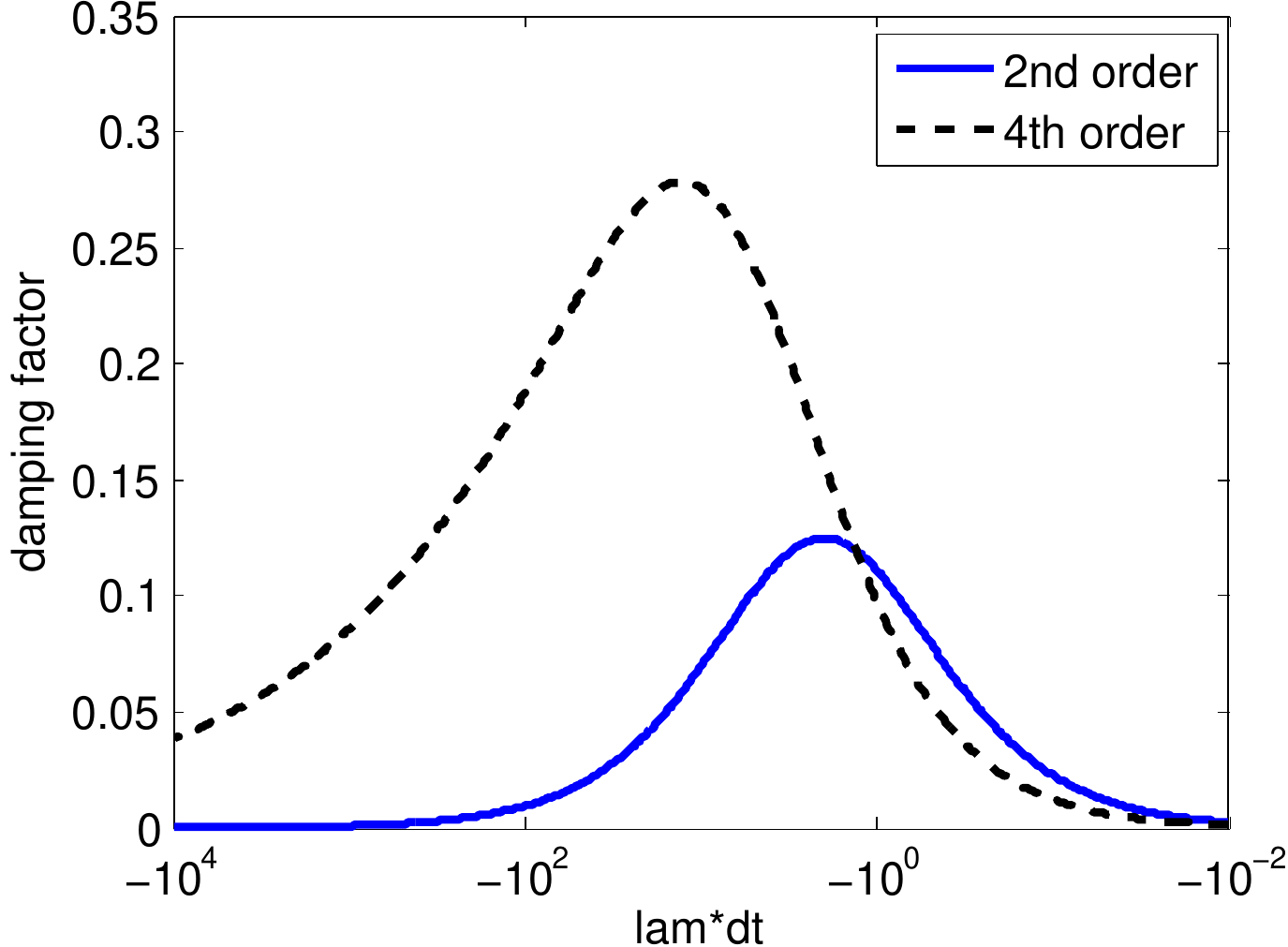}
\caption{Damping factors for second- and fourth-order SDC methods
used in Section \ref{sect:numerics}.}
\end{figure}

Figure \ref{fig:damping} shows the maximum magnitude of the 
eigenvalues for both second- and fourth-order cases as a function of 
$\lamDt$.  Clearly in  both limits
$\lamDt \rightarrow 0$ and $\lamDt \rightarrow -\infty$, 
the damping factor goes to zero. This behavior differs from the 
traditional analysis of point relaxation for elliptic equations
like the Laplace equation where high-frequency eigen-components are 
damped more rapidly than low frequency components for the classical
relaxation schemes.  This feature is the reason why multigrid methods
provide a tremendous speedup compared to relaxation alone. 
In the context of SDC sweeps, Figure \ref{fig:damping} demonstrates
that although MLSDC with temporal coarsening may provide some 
increase in efficiency compared to higher-order SDC, SDC alone
will still converge rapidly in the stiff limit for this problem.

For any linear problem, one can represent one iteration of PFASST 
or IPFASST as a matrix vector multiplication on a vector composed of the solution at
each SDC node within each parallel time step.  Examining the eigenvalues
of this matrix would then give an indication of how PFASST 
iterations would converge for a given choice of parameters.  One could
then examine,  for a given linear PDE,  different numbers and types of 
nodes; types of sweeps and quadrature rules; number of levels,
refinement factors in space and time, and  type of interpolation and restriction
between levels; and the size and  number of parallel time-steps.  In the 
case of IPFASST, the type of relaxation, the order of 
interpolation and restriction, and the  number of V-cycles per implicit-solve
could also be varied\footnote{It is also worth noting that the use of non-uniform sub-steps
in SDC removes the equivalence between Fourier and eigenmode decomposition
in the time direction.}.  
A paper presenting a systematic study of the linear convergence of PFASST and IPFASST
is in preparation.

%% file: numerics.tex
\section{Numerical Results}\label{sect:numerics}

In this section, numerical results illustrating some
features of the  IPFASST method are presented on the
model problem given by the linear heat equation.  First,
results analyzing the convergence behavior of IPFASST
are reported.  Strong scaling timing results for three-dimensional 
problems using space-time parallelism are then presented
in Section \ref{sect:strong}.

Despite the popularity of the heat equation as a test case for space-time
parallel methods
(e.g. \cite{HortonKnirsch1992,HortonVandewalle1995,VandewalleHorton1995,Gander1998,AmodioBrugnano2009,FriedhoffEtAl2013}), 
there are several reasons why it is not a particularly good
test case for time parallelization.  First, when considering time-parallel
methods, one would like to show how the method scales as the number of parallel
time steps grows.  Hence for many time steps, one must either choose a 
relatively long interval of integration or a relatively small time step.
If a long time interval is chosen, the solution will decay to a value
close to zero which complicates any discussion of accuracy and convergence.
If on the other hand, very small time steps are chosen, then the temporal
accuracy of the method will likely far exceed the spatial accuracy, which
brings into question why such a small time step is being considered.
Finally, it is important to note that 
the performance of parallel-in-time methods applied to the the linear heat
equation may not be indicative of performance on other problems that
are not strictly parabolic, have dynamic spatial features, or have nontrivial
boundary conditions.
Despite these drawbacks, we will investigate the performance of IPFASST 
on this test case.  The main  motivation for doing so  is 
to provide a comparison of the PFASST method with other published
space-time parallel methods using this example.

\subsection{One-Dimensional Convergence Studies} \label{sect:oneD}

In this section we consider the accuracy and convergence
of IPFASST, including the dependence on the number of V-cycles.
For reasons addressed below, we use here the simple model
problem consisting of the one-dimensional heat equation
with homogeneous Dirichlet boundary conditions. 
The equation in a general form prescribed
on the space-time domain  $[0,L] \times [0,T]$  is given by
\begin{align*}
  u_t   &=  \nu u_{xx} \\
  u(x,0) &= u_0(x) \\
  u(0,t)  &= u(L,t) = 0.
\end{align*}
Choosing the initial conditions
\begin{displaymath}
  u_0(x)  =  \sin(k \pi x/L),
\end{displaymath}
gives the exact solution
\begin{displaymath}
  u(x,t) =  e^{-\nu  (k \pi/L)^2 t} u_0(x).
\end{displaymath}
Using the method of lines and a second-order finite-difference approximation of the spatial derivatives  gives the linear system of ODEs
\begin{displaymath} \label{heat:ODE}
  u'_i(t) =  \nu \frac{u_{i-1}-2 u_i + u_{i+1}}{\Delta x^2},
\end{displaymath}
where $u_i(t) \approx u(i \Delta x, t)$ for $i=1 \ldots N-1$, $\Delta x = L/N$, and $u_0=u_N=0$.
The exact solution of the systems of ODEs given initial conditions
\begin{displaymath}
  u_i(0)  =  \sin(k \pi x_i/L),
\end{displaymath}
is
\begin{displaymath}
  u_i(t) =  e^{-\nu d(k) t} u_0(x_i),
\end{displaymath}
where
\begin{displaymath}
  d(k) =  \frac{-2 + 2\cos(k \pi \Dx/L)}{\Dx^2}.
\end{displaymath}
Note that for $k$ even moderately large, the solution decays very rapidly to zero.
In the following one-dimensional examples, we choose $k=L=\nu=T=1$, which
means the solution decays to approximately $5\times 10^{-5}$ during the time interval.

One of the convenient features of SDC methods is that it
is simple to construct higher-order methods simply by increasing the
number of quadrature nodes being used.  Even restricting the
discussion to second-order spatial discretizations, using a
first-order or second-order temporal integration method is very
inefficient.  As a simple demonstration of this, consider the
following numerical experiment.  We apply SDC methods of different
orders to \eqref{heat:ODE} with $k=\nu=L=T=1$. Specifically we 
consider SDC methods of order 1,2,4, and 8 where the number of SDC
nodes per time step corresponds to the formal order.
We chose $\Dx=1/128$
and compute the $L_\infty$ error of the solution at the final time
$T=1$ compared to the exact solution of the PDE and the exact solution
of the discrete ODE for various values of $\Dt = 1/\Nt$ for $\Nt=2^p$ with
$p$ ranging from 2 to 12.  In addition, we compute the residual at the final
time step as in \eqref{eq:resid}.

\begin{figure}[t]
\centering
 \includegraphics[width=0.95\textwidth]{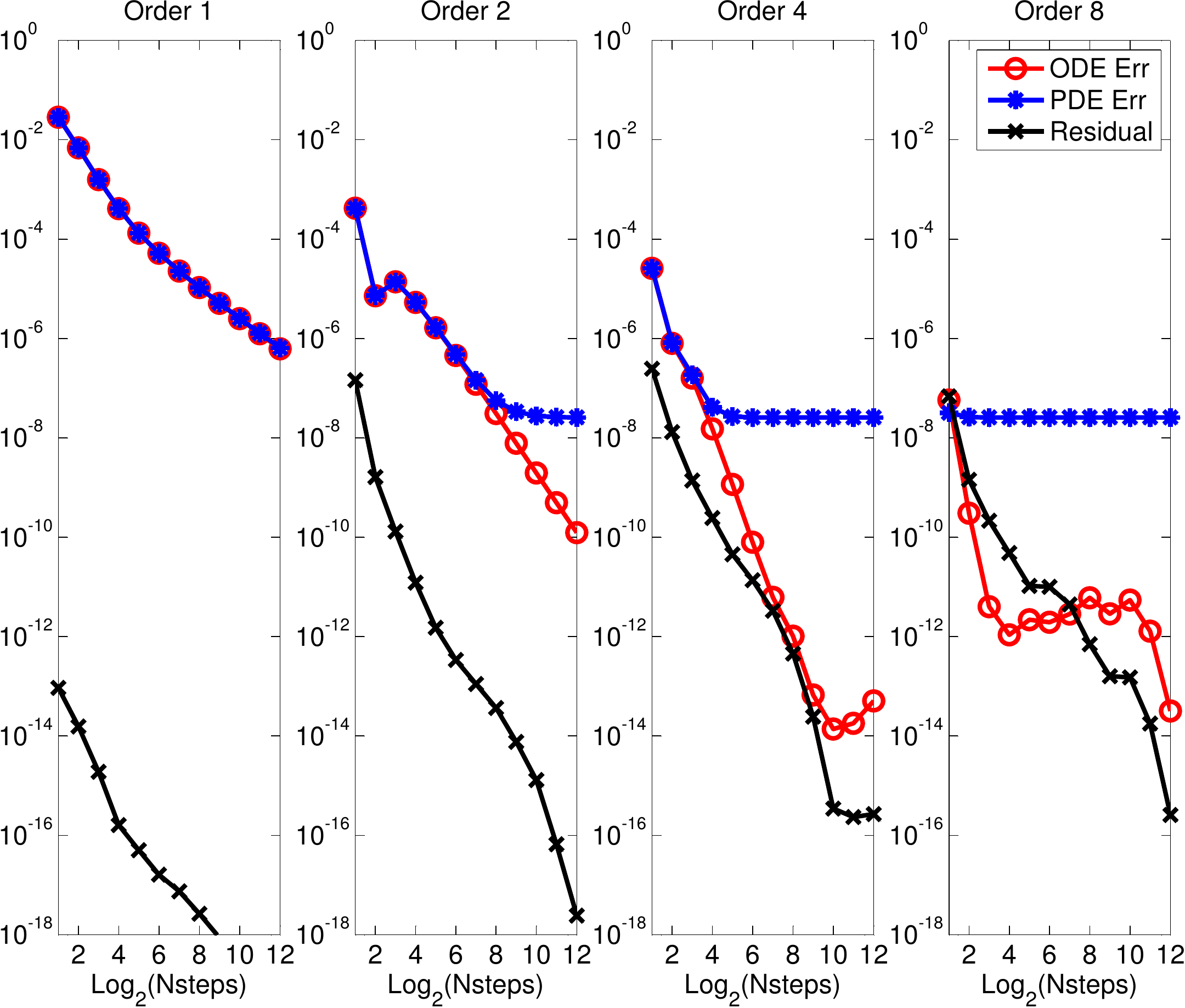}
\caption{Convergence in error and residual for serial SDC methods of different orders
for the 1-D heat equation.} \label{fig:Conv3}
\end{figure}

The results are displayed in Figure \ref{fig:Conv3}.  Several comments can be made.
Most obviously, as the order of the time integration increases, the number of 
time steps required to reduce the PDE error to the minimum possible given
the spatial resolution goes down dramatically.  
For the eighth-order method, two
time steps of size 0.5 yields a smaller error than 128 time steps of the second-order
method.  
Of course, higher-order methods
require more work per time step due to the increase in substeps on the SDC nodes, but
as shown in Section \ref{sect:strong} higher-order discretizations in
space and time lead to much reduced computational times for problems in three
dimensions.  Note also that the residual and error in the discretized ODE continue to
converge to zero well past the minimum error due to the spatial discretization.  
This implies that a judicious choice of residual tolerance is needed for iterative
approaches like IPFASST to avoid wasted iterations.

In the PFASST algorithm, parallel efficiency depends on the availability of
a hierarchy of space-time discretizations as in MLSDC, where coarsening in the temporal
direction is achieved by reducing the number of nodes used in the underlying
SDC method on each processor.  Obviously for a second-order temporal discretization,
only one coarsening step in time is possible leading to a coarse level based on
the backward-Euler method. In the spatial directions, one way to reduce the work on 
PFASST levels is to reduce the order of the spatial approximation on coarser 
levels, however, this is only possible when higher-order spatial approximation is 
being used on the finest level.  Hence the options for coarsening in PFASST
when restricted to second-order discretizations are more limited than for 
higher-order methods.

\subsubsection{Comparison of IPFASST With Different Number of V-cycles}

In this numerical experiment, we compare the convergence of IPFASST
with different numbers of prescribed V-cycles for each approximate
linear solve at each sub-step. A standard geometric multigrid method
is used for these tests with linear interpolation,
full-weighting restriction,
and two pre- and two post-smoothing sweeps using a Jacobi smoother.
IPFASST will only be more efficient
than PFASST if the reduced cost of fewer V-cycles is not offset by
an increase in the number of iterations required for convergence.

\begin{figure}[t]
	\centering
        \includegraphics[width=0.95\textwidth]{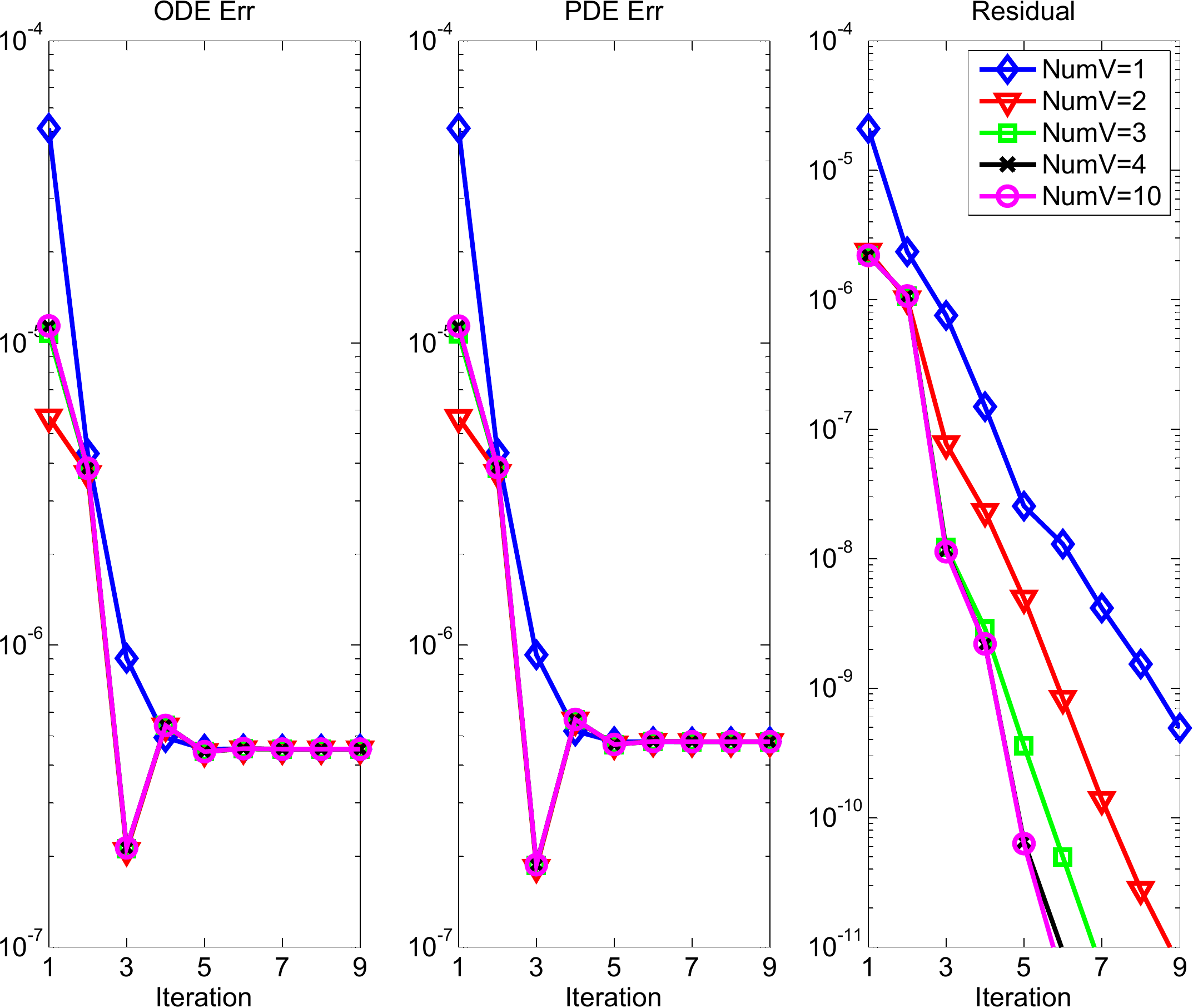}
	\caption{Comparison of the convergence of IPFASST using different numbers of V-cycles per solve.}\label{fig:vcycles}
\end{figure}

We use here $\Nx=\Nt=128$ with 128 parallel time steps and 
the same parameters as in the previous examples, namely $k=\nu=L=T=1$.
Figure~\ref{fig:vcycles} shows the convergence of IPFASST in terms of ODE
error, PDE error and residual versus the number of V-cycles used
for the approximate implicit solve in every sub-step.  The residual
decays faster as the number of V-cycles is increased from $1$ to $3$,
but using more than $3$ V-cycles does not yield further improvement.
In contrast, both the ODE and the PDE error are less effected by the
number of V-cycles.  Using  two V-cycles per solve gives virtually
the same behavior as using 10. We stress again that these results
for the linear heat equation may not hold for PDEs of different
mathematical type (see \cite{SpeckEtAl2014_DDM2013} for some
preliminary results on advection diffusion equations).

\subsubsection{Convergence of IPFASST With Weak Scaling}
An important motivation for the study of time-parallel methods is the
``trap of weak scaling" given current supercomputer evolution where
cores counts are increasing, but processor speed is not.  For a given
application, if the spatial resolution is increased as core counts
increase, the cost of spatial operators will remain close to constant
assuming good weak scaling of the spatial operations.  However, the
time step size will eventually need to also be reduced, either because of
stability constraints and/or to match the increased spatial accuracy
(depending on whether explicit or implicit methods are used).  This
means that more time steps are necessary for the same simulation time,
and therefore, the total runtime will increase even with perfect
spatial scaling unless some parallelism in the time direction is
employed.  For time-parallel methods to be effective, it is necessary
that the performance does not deteriorate as resolution is refined.
Therefore, in this test we consider how the convergence of the IPFASST
iterations scales with the problem size and number of parallel time
steps.

Three runs are performed with $N_x=N_t=32$, $64$, and $128$
again choosing $k=\nu=L=T=1$.  The
number of parallel time-steps in IPFASST is equal to the number of time
steps $N_t$.  We use three levels with coarsening by a factor two in space
and second-order finite-difference approximations on all levels.  On
the coarsest level, the collocation rule corresponds to first-order
backward Euler.  On the middle and finest level, we use $3$ uniform
nodes corresponding to a second-order collocation rule.  Based on
the previous experiments, implicit solves
are approximated by two V-cycles with two pre- and post-smoothing
steps with a Gauss-Seidel smoother.  

\begin{figure}[t]
	\centering
    	\includegraphics[width=0.95\textwidth]{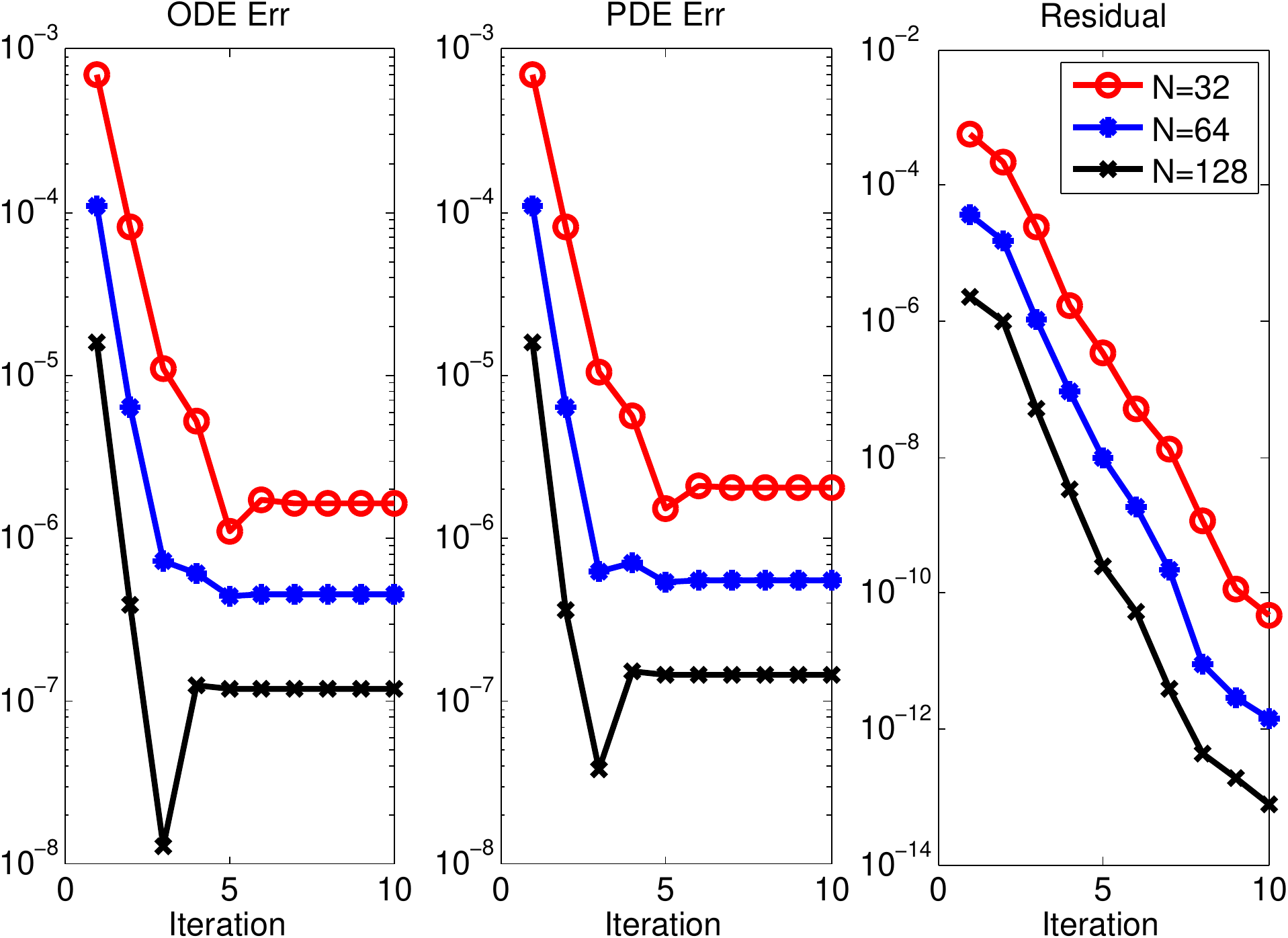}
	\caption{Convergence of IPFASST in terms of ODE error, PDE error, and residual in the final time step for increasing resolution in both space and time.}\label{fig:optimal}
\end{figure}

Figure~\ref{fig:optimal} shows how IPFASST converges for the different
resolutions in terms of the ODE error (left), the PDE error (middle)
and the residual (right).  Errors are reported for the last time step.
As the resolution 
increases in space or time, 
the error level decreases up to some minimum level.
Depending on which discretization error is
dominant, this level is either the discretization error of the
collocation rule or of the spatial discretization, cf. the discussion
in \cite[Sect. 3.2]{SpeckEtAl2014_BIT}.  This level is reached at
iteration $5$ for $N=32$ and at iteration $3$ for $N=64$ and $N=128$.
Therefore, increasing the resolution and with it the number of
parallel time steps does not increase the number of iterations
required by IPFASST.  Note that this does not cover the case where the
number of concurrently computed steps is increased while
the time step
size is kept constant in order to compute over a longer time interval.

\subsection{Three-dimensional Strong Scaling Studies} \label{sect:strong}
To compare the computational cost of  IPFASST with serial methods 
we consider the 3D heat equation
\begin{equation}
	\label{eq:heat3d}
	u_t(\tvect{x},t) = \nu \Delta u(\tvect{x}, t), \quad \tvect{x} \in \Omega, \quad 0 \leq t \leq T
\end{equation}
on the unit cube $\Omega = [ 0, 1 ]^{3}$ with $T=1.0$, initial conditions
\begin{equation}
	u(\tvect{x},0) = \sin(\pi x) \sin(\pi y) \sin(\pi z),
\end{equation}
and homogeneous Dirichlet boundary conditions.  We choose $\nu =
\frac{1}{3}$ so that the solution decays at the same rate as the
previous one-dimensional tests.  The Laplacian is discretized with
either a second- or fourth-order finite difference stencil, and either
ISDC or IPFASST is used to solve the resulting initial value problem.
In all cases an implicit Euler sub-stepping is used for the SDC sweeps,
and  PMG \cite{bolten:2014} is used for parallel multigrid V-cycles.
Simulations are run on the IBM BlueGene/Q installation JUQUEEN at
J\"ulich Supercomputing Centre.  Timing results are displayed
in Figure~\ref{fig:ipfasst_timing}, and the specifications of the 
different runs are provided below.

\subsubsection{Problem and Method Specifications}
ISDC, IMLSDC and IPFASST are run in three different configurations,
which all give approximately identical errors measured in the
infinity norm against the analytical solution of~\eqref{eq:heat3d}.
The tolerance for the residual of ISDC and IPFASST is set to
$10^{-9}$, which in all three configurations results in an error of
about $1.5 \times 10^{-7}$.  PMG  uses a tolerance of $10^{-12}$ 
whenever a full solve is performed,
and a stalling criterion that stops the iteration if the new residual
is not smaller then $75 \%$ of the previous one.  Two levels are used
in IMLSDC and IPFASST, with  fourth-order spatial interpolation and
point-wise restriction in both space and time.

In all runs, approximate implicit solves consisting of two PMG V-cycles are 
used in the SDC sweeps. PMG V-cycles use two pre- and two post-smoothing 
sweeps consisting of JOR red-black smoothers.  Linear interpolation and 
full-weighting restriction are used in the V-cycles.

\paragraph{Second-order spatial and second-order temporal
  discretization} (Marked by circles in
Figure~\ref{fig:ipfasst_timing}). For these runs, ISDC and the fine
level of IMLSDC/IPFASST use two quadrature nodes corresponding to the
midpoint of the time step and the right-hand value $t_{n+1}$. 
Therefore the method converges to a second-order collocation scheme. 
The coarse levels in IMLSDC and IPFASST correspond to backward Euler. The
spatial mesh has $N = 128^{3}$ points and $128$ time steps are
performed. The coarse level uses a $64^3$ point spatial mesh. For
ISDC, runs are performed with the number of cores used by PMG varying
between $16$ and $32,768$.  For IMLSDC and IPFASST, the number of
cores for PMG is fixed to $4,096$ and the number of parallel time-ranks varied
between $1$ (for IMLSDC) up to $32$, for a total of $4,096 \times 32 =
131,072$ cores.  On average, ISDC requires about $3.5$ iterations and
IMLSDC about $3.7$ for convergence.  The last time-rank in IPFASST requires between
about $3.4$ (for two parallel time-ranks) and $4.0$ iterations (for
$32$ time-ranks).
\paragraph{Second-order spatial and fourth-order temporal discretization} 
(Marked by diamonds in Figure~\ref{fig:ipfasst_timing}). Here, ISDC
and the fine level of IPFASST use 4 quadrature nodes, corresponding to
fourth-order collocation.  The coarse level of IPFASST is the same as
the second-order runs above.  Only $24$ time steps are required to
achieve the same overall error because of the higher-order temporal
discretization. The spatial mesh remains as described above.  As
before, scaling of ISDC is measured using $16$ up to $32,768$ cores
for PMG while the number of cores for PMG in IMLSDC and IPFASST is
fixed at $4,096$.  Because here only $24$ time steps have to be
performed, the number of time-ranks is varied only from $1$ up to
$24$, for a total maximum number of cores of $24 \times 4,096 =
98,304$ cores.  IMLSDC needs an average of $3.8$ iterations, and the
last time-rank in IPFASST between $4.0$ (for two time-ranks) up to
$6.0$ (for $24$ time-ranks).
\paragraph{Fourth-order spatial and fourth-order temporal  discretization} 
(Marked by squares in Figure~\ref{fig:ipfasst_timing}).  As above,
$24$ time-steps are used with a fourth-order temporal discretization.
$N=32^{3}$ spatial points suffice to maintain the same error as above. On the
coarse level, IPFASST uses only $N=16^3$ points and a second-order
finite-difference stencil. Here, because much fewer spatial degrees of
freedom are needed, scaling of ISDC is measured only using $16$ to
$512$ cores in PMG.  IMLSDC and IPFASST use $64$ cores for PMG and, as
before, up to $24$ time-ranks.  The maximum total number of cores here
is therefore only $24 \times 64 = 1,536$.

\subsubsection{Results}
\begin{figure}[t]
\label{fig:Qruns}
\centering
    \includegraphics[width=0.95\textwidth]{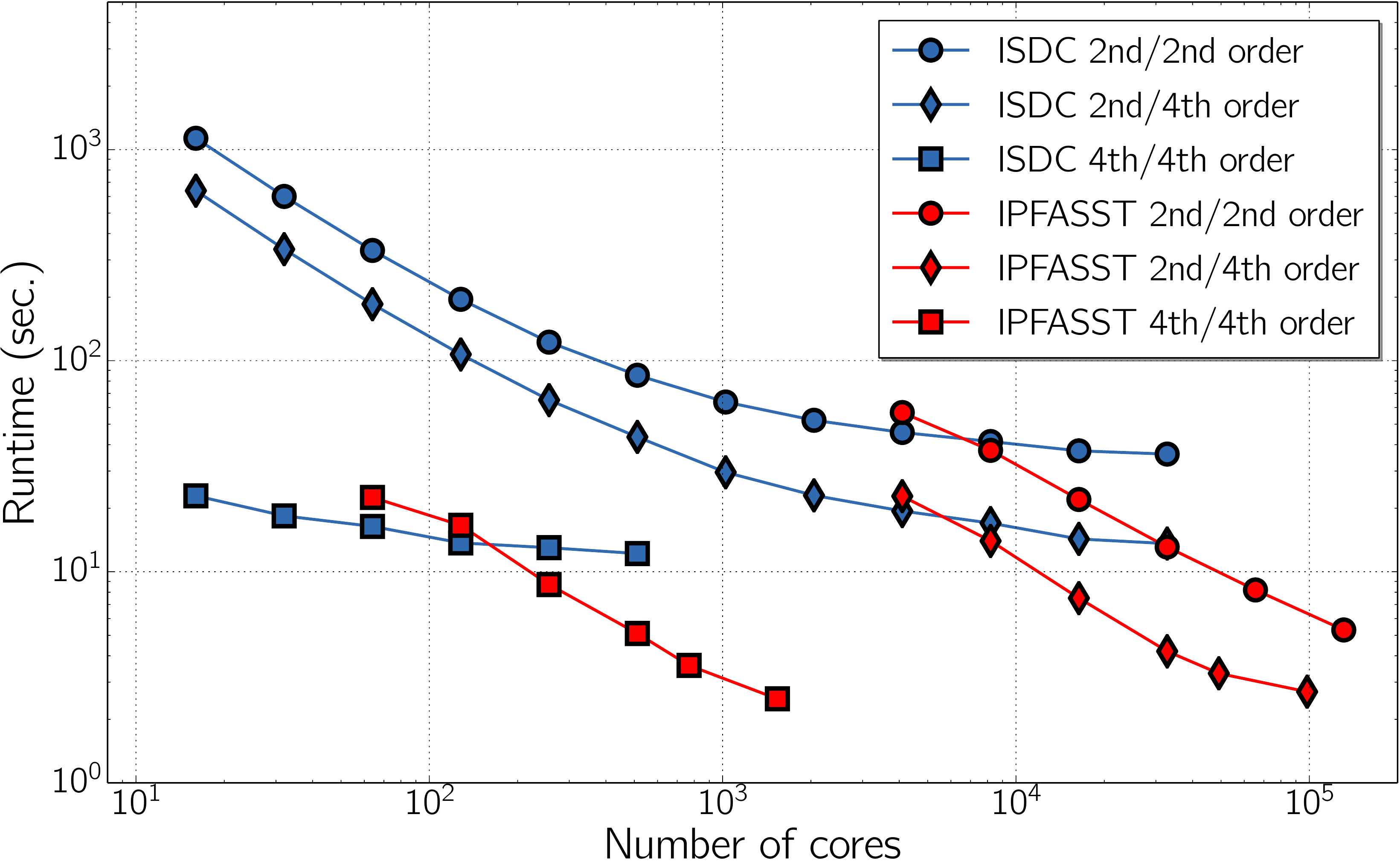}
\caption{Run times for the 3D heat equation depending on the total number of cores used
for different configurations of ISDC and IPFASST. The first entry for each configuration 
of IPFASST corresponds to a run with only a single-time rank, that is IMLSDC. 
All setups result in an error of about $1.5 \times 10^{-7}$, making a 
comparison of their runtimes meaningful.}\label{fig:ipfasst_timing}
\end{figure}
Figure~\ref{fig:ipfasst_timing} displays the run times for the different
configurations of ISDC and IPFASST laid out above.  In addition,
Table~\ref{tab:ipfasst} lists parallel speedup and efficiency 
(speedup divided by number of processors) provided by
IPFASST compared to the ISDC counterpart with the same number of spatial
processors for PMG.  Note that the first marker for each IPFASST line
corresponds to a run with only a single time-rank, which is  IMLSDC.
For each of the three setups, IMLSDC by itself is somewhat  slower
than ISDC using the same number of cores for PMG due to a slight increase
in the number of iterations required by IMLSDC and the overhead of 
coarsening and interpolation between SDC levels.  This is in contrast
to the SDC/MLSDC method using full solves on the fine level as studied
in~\cite{SpeckEtAl2014_BIT}, where MLSDC could significantly reduce
run times compared to SDC.   This means that the most efficient parallel
variant for this test case (namely IPFASST) is not a direct parallelization
of the most efficient serial variant (namely ISDC).  Hence while
the use of incomplete solves reduces the computational cost for both
ISDC and IPFASST, it actually decreases the parallel efficiency of IPFASST
because we have to compare to ISDC.  We should also note that ISDC
is not necessarily the most efficient serial method for this problem,
but we use it as a comparison to the parallel performance of IPFASST
based on ISDC sweeps.

As already observed in~\cite{SpeckEtAl2014_parco} for the PFASST
method, it is more efficient to  begin space-time parallelism with
IPFASST using fewer spatial cores for PMG as for when the 
parallel speedup saturates.  At the limit of spatial scaling,
the spatial coarsening within PFASST does 
not  make the coarse levels sweeps much cheaper than the fine level and
therefore reduces the parallel efficiency of PFASST or IPFASST.

While the second/second-order and the second-/fourth-order versions
of ISDC and IPFASST scale approximately equally well, the higher-order
methods result in shorter run times.  For fully fourth-order ISDC,
spatial scaling is of course significantly worse than for the 
second-order spatial methods 
because the size of the problem is drastically smaller. However,
no matter how many processors are used for the second-order version,
the fourth-order ISDC is always significantly faster.  The same is
true for IPFASST. While all three configurations of IPFASST scale
approximately equally well, the fourth-order version requires
significantly fewer cores to achieve the same runtime as the
second-order version.  The smallest time-to-solution is about $2.5$
seconds provided by fourth-order IPFASST using a total of $1,536$
cores.  While the mixed second/ fourth-order version is only slightly
slower at $2.7$ seconds, it requires $98,304$ 
cores to achieve this runtime.  We stress again that these results may not
translate to other problems, particularly where higher-order methods
are not appropriate due to a lack of smoothness in the solutions.
\begin{table}[t]
	\centering
	\footnotesize
	\begin{tabular}{|c|c|c|}\hline
		\multicolumn{3}{|c|}{\2nd / \2nd order } \\ \hline
		$N_{p}$ & Speedup & Efficiency \\ \hline
		1   & 0.8 & -- \\
		2   & 1.2 & 60.0 \% \\
		4   & 2.1 & 53.5 \% \\
		8   & 3.5 & 43.5 \% \\
		16 & 5.6 & 35.0 \% \\
		32 & 8.6 & 26.9 \% \\ \hline
	\end{tabular}
	\begin{tabular}{|c|c|c|}\hline
		\multicolumn{3}{|c|}{\2nd / \4th order } \\ \hline
		$N_{p}$  & Speedup & Efficiency \\ \hline
		1   & 0.9 & -- \\
		2   & 1.4  & 70.0 \% \\
		4   & 2.6  & 65.0 \%\\
		8   & 4.6  & 57.5 \% \\
		12 & 5.9  & 49.2 \% \\
		24 & 7.2  & 30.0 \% \\ \hline
	\end{tabular}
	\begin{tabular}{|c|c|c|}\hline
		\multicolumn{3}{|c|}{\4th / \4th order } \\ \hline
		$N_{p}$  & Speedup & Efficiency \\ \hline
		1   & 0.7 & -- \\
		2   & 1.0 & --  \\
		4   & 1.9 & 47.5 \% \\
		8   & 3.2 & 40.0 \% \\
		12 & 4.6 & 38.3 \% \\
		24 & 6.6 & 27.5 \% \\ \hline		
	\end{tabular}		
	\caption{Additional speedup provided by IPFASST compared to the ISDC run with the same number of spatial processors  for PMG. $N_{p}$ indicates the number of parallel time steps.}\label{tab:ipfasst}
\end{table}

%% file: discussion.tex
\section{Discussion and Outlook}

The two main goals for 
this paper are (1)
to investigate combining spatial and temporal iterative strategies for
space-time parallelization and (2) to provide some data on the
performance of the resulting IPFASST method applied to the linear heat equation
for comparison with other published results. The IPFASST algorithm
introduced here interweaves the iterations of a spatial multigrid
solver (PMG) with the temporal iterations in the SDC methods to achieve
space-time parallelism, where the resulting hybrid space-time
iterations are constructed by considering the spatial and temporal
dimensions independently. Numerical results suggest that reducing
the number of V-cycles for implicit spatial solves in the IPFASST
can be done without a significant impact on the convergence of
the time-parallel iterations (see e.g. Fig. 5). 
The scaling results shown in Fig. 7 
demonstrate that incorporating temporal parallelization as in 
IPFASST can extend strong scaling and further reduce
the time-to-solution when spatial parallelism is close to saturation, 
assuming  more resources are available.

One interesting result suggested by the numerical experiments is 
that although MLSDC might be more efficient than SDC for the heat
equation, it seems that ISDC is more efficient than IMLSDC. The
difference between ISDC and IMLSDC is essentially only the order in
which SDC relaxation sweeps and multigrid V-cycles are performed on
various levels, and hence it is likely that even more efficient
variations using a more general ordering of space-time relaxations
could be found. A careful analysis of different variations of the
methods presented here, and the extension of the linear analysis in
Sect.~3 to PFASST and IPFASST is currently underway.  As mentioned,
given the number of different options already possible, finding
optimal configurations may be difficult, and more importantly these
optimal configurations are probably problem-dependent.

%% file: ack.tex
\section{Acknowledgments}

Authors Minion and Emmett were supported by the Applied Mathematics
Program of the DOE Office of Advanced Scientific Computing Research
under the U.S. Department of Energy under contract DE-AC02-05CH11231.
Speck and Ruprecht acknowledge support by the Swiss National Science Foundation (SNSF) under the lead agency agreement through the project ``ExaSolvers'' (SNF-145271) within the Priority Programme 1648 ``Software for Exascale Computing'' (SPPEXA) of the Deutsche Forschungsgemeinschaft. 
Computing time on JUQUEEN at the J\"ulich Supercomputing
Centre (JSC) was provided by project HWU12.

%% file: paper.bbl
\begin{thebibliography}{10}

\bibitem{AmodioBrugnano2009}
{\sc P.~Amodio and L.~Brugnano}, {\em Parallel solution in time of odes: some
  achievements and perspectives}, Applied Numerical Mathematics, 59 (2009),
  pp.~424 -- 435.

\bibitem{bohmerHemkerStetter:1984}
{\sc K.~B\"ohmer, P.~Hemker, and H.~J. Stetter}, {\em The defect correction
  approach}, in Defect Correction Methods. Theory and Applications, K.~B\"ohmer
  and H.~J. Stetter, eds., Springer-Verlag, 1984, pp.~1--32.

\bibitem{bolten:2014}
{\sc M.~Bolten}, {\em Evaluation of a multigrid solver for 3-level {Toeplitz}
  and circulant matrices on {Blue Gene/Q}}, in NIC Symposium 2014, K.~Binder,
  G.~M{\"u}nster, and M.~Kremer, eds., John von Neumann Institute for
  Computing, 2014, pp.~345--352.
\newblock (to appear).

\bibitem{BourliouxEtAl2003}
{\sc Anne Bourlioux, Anita~T. Layton, and Michael~L. Minion}, {\em High-order
  multi-implicit spectral deferred correction methods for problems of reactive
  flow}, Journal of Computational Physics, 189 (2003), pp.~651 -- 675.

\bibitem{BouzarthMinion2010}
{\sc Elizabeth~L. Bouzarth and Michael~L. Minion}, {\em A multirate time
  integrator for regularized {S}tokeslets}, Journal of Computational Physics,
  229 (2010), pp.~4208 -- 4224.

\bibitem{multigridtutorial}
{\sc William~L. Briggs, Van~Emden Henson, and Steve~F. Mc{C}ormick}, {\em A
  Multigrid Tutorial}, SIAM, 2000.

\bibitem{Burrage1997}
{\sc Kevin Burrage}, {\em {Parallel methods for {ODE}s}}, Advances in
  Computational Mathematics, 7 (1997), pp.~1--3.

\bibitem{DuttEtAl2000}
{\sc Alok Dutt, Leslie Greengard, and Vladimir Rokhlin}, {\em Spectral deferred
  correction methods for ordinary differential equations}, BIT Numerical
  Mathematics, 40 (2000), pp.~241--266.

\bibitem{EmmettMinion2012}
{\sc M.~Emmett and M.~L. Minion}, {\em Toward an efficient parallel in time
  method for partial differential equations}, Communications in Applied
  Mathematics and Computational Science, 7 (2012), pp.~105--132.

\bibitem{EmmettMinion2014_DDM}
{\sc Matthew Emmett and Michael~L. Minion}, {\em Efficient implementation of a
  multi-level parallel in time algorithm}, in Domain Decomposition Methods in
  Science and Engineering XXI, vol.~98 of Lecture Notes in Computational
  Science and Engineering, Springer International Publishing, 2014,
  pp.~359--366.

\bibitem{FalgoutEtAl2014_MGRIT}
{\sc R.~D. Falgout, S.~Friedhoff, Tz.~V. Kolev, S.~P. MacLachlan, and J.~B.
  Schroder}, {\em Parallel time integration with multigrid}, SIAM Journal on
  Scientific Computing, 36 (2014), pp.~C635--C661.

\bibitem{FriedhoffEtAl2013}
{\sc S~Friedhoff, R~D Falgout, T~V Kolev, S~MacLachlan, and J~B Schroder}, {\em
  A multigrid-in-time algorithm for solving evolution equations in parallel},
  in Presented at: Sixteenth Copper Mountain Conference on Multigrid Methods,
  Copper Mountain, CO, United States, Mar 17 - Mar 22, 2013, 2013.

\bibitem{Gander1998}
{\sc Martin~J. Gander and Andrew~M. Stuart}, {\em {Space-Time Continuous
  Analysis of Waveform Relaxation for the Heat Equation}}, SIAM Journal on
  Scientific Computing, 19 (1998), pp.~2014--2031.

\bibitem{Geist2010}
{\sc Al~Geist}, {\em {Paving the roadmap to exascale}}, SciDAC Review, 16
  (2010), pp.~53--59.

\bibitem{Hackbusch1984}
{\sc W.~Hackbusch}, {\em Parabolic multi-grid methods}, Computing Methods in
  Applied Sciences and Engineering, VI,  (1984), pp.~189--197.

\bibitem{HagstromZhou2006}
{\sc Thomas Hagstrom and Ruhai Zhou}, {\em On the spectral deferred correction
  of splitting methods for initial value problems}, Communications in Applied
  Mathematics and Computational Science, 1 (2006), pp.~169--205.

\bibitem{HortonKnirsch1992}
{\sc G.~Horton and R.~Knirsch}, {\em A time-parallel multigrid-extrapolation
  method for parabolic partial differential equations}, Parallel Computing, 18
  (1992), pp.~21 -- 29.

\bibitem{HortonVandewalle1995}
{\sc G.~Horton and S.~Vandewalle}, {\em A space-time multigrid method for
  parabolic partial differential equations}, SIAM Journal on Scientific
  Computing, 16 (1995), pp.~848--864.

\bibitem{HuangEtAl2006}
{\sc Jingfang Huang, Jun Jia, and Michael Minion}, {\em Accelerating the
  convergence of spectral deferred correction methods}, Journal of
  Computational Physics, 214 (2006), pp.~633 -- 656.

\bibitem{LaytonMinion2004}
{\sc Anita~T. Layton and Michael~L. Minion}, {\em Conservative multi-implicit
  spectral deferred correction methods for reacting gas dynamics}, Journal of
  Computational Physics, 194 (2004), pp.~697 -- 715.

\bibitem{LaytonMinion2005}
{\sc Anita~T. Layton and Michael~L. Minion}, {\em Implications of the choice of
  quadrature nodes for picard integral deferred corrections methods for
  ordinary differential equations}, BIT Numerical Mathematics, 45 (2005),
  pp.~341--373.

\bibitem{LionsEtAl2001}
{\sc J.-L. Lions, Y.~Maday, and G.~Turinici}, {\em A ''parareal'' in time
  discretization of {P}{D}{E}'s}, Comptes Rendus de l'Acad\'{e}mie des Sciences
  - Series I - Mathematics, 332 (2001), pp.~661--668.

\bibitem{Minion2003}
{\sc Michael~L. Minion}, {\em Semi-implicit spectral deferred correction
  methods for ordinary differential equations}, Communications in Mathematical
  Sciences, 1 (2003), pp.~471--500.

\bibitem{Minion2010}
{\sc M.~L. Minion}, {\em A hybrid parareal spectral deferred corrections
  method}, Communications in Applied Mathematics and Computational Science, 5
  (2010), pp.~265--301.

\bibitem{MinionEtAl2008}
{\sc Michael~L Minion and Sarah~A Williams}, {\em Parareal and spectral
  deferred corrections}, in AIP Conference Proceedings, vol.~1048, 2008,
  p.~388.

\bibitem{Nievergelt1964}
{\sc J.~Nievergelt}, {\em Parallel methods for integrating ordinary
  differential equations}, Commun. ACM, 7 (1964), pp.~731--733.

\bibitem{Pereyra1966}
{\sc Victor Pereyra}, {\em {On improving an approximate solution of a
  functional equation by deferred corrections}}, Numerische Mathematik, 8
  (1966), pp.~376--391.

\bibitem{RuprechtEtAl2013_SC}
{\sc D.~Ruprecht, R.~Speck, M.~Emmett, M.~Bolten, and R.~Krause}, {\em Poster:
  Extreme-scale space-time parallelism}, in Proceedings of the 2013 Conference
  on High Performance Computing Networking, Storage and Analysis Companion, SC
  '13 Companion, 2013.

\bibitem{SpeckEtAl2014_parco}
{\sc R.~Speck, D.~Ruprecht, M.~Emmett, M.~Bolten, and R.~Krause}, {\em A
  space-time parallel solver for the three-dimensional heat equation}, in
  Parallel Computing: Accelerating Computational Science and Engineering (CSE),
  vol.~25 of Advances in Parallel Computing, IOS Press, 2014, pp.~263 -- 272.

\bibitem{SpeckEtAl2014_BIT}
{\sc R.~Speck, D.~Ruprecht, M.~Emmett, M.~Minion, M.~Bolten, and R.~Krause},
  {\em A multi-level spectral deferred correction method}, {BIT} Numerical
  Mathematics,  (2014).

\bibitem{SpeckEtAl2012}
{\sc R.~Speck, D.~Ruprecht, R.~Krause, M.~Emmett, M.~Minion, M.~Winkel, and
  P.~Gibbon}, {\em A massively space-time parallel {N}-body solver}, in
  Proceedings of the International Conference on High Performance Computing,
  Networking, Storage and Analysis, SC '12, Los Alamitos, CA, USA, 2012, IEEE
  Computer Society Press, pp.~92:1--92:11.

\bibitem{SpeckEtAl2014_DDM2013}
{\sc R.~Speck, D.~Ruprecht, M.~Minion, M.~Emmett, and R.~Krause}, {\em Inexact
  spectral deferred corrections}, in Domain Decomposition Methods in Science
  and Engineering {XXII}, Lecture Notes in Computational Science and
  Engineering, 2014.
\newblock Accepted.

\bibitem{stetter:1974}
{\sc H.~J. Stetter}, {\em Economical global error estimation}, in Stiff
  Differential Systems, R.~A. Willoughby, ed., 1974, pp.~245--258.

\bibitem{VandewalleHorton1995}
{\sc S.~Vandewalle and G.~Horton}, {\em Fourier mode analysis of the multigrid
  waveform relaxation and time-parallel multigrid methods}, Computing, 54
  (1995), pp.~317--330.

\bibitem{Weiser2013}
{\sc Martin Weiser}, {\em Faster {SDC} convergence on non-equidistant grids
  with {DIRK} sweeps}.
\newblock {ZIB} Report 13--30, 2013.

\bibitem{Zadunaisky1966a}
{\sc P.~E. Zadunaisky}, {\em A method for the estimation of errors propagated
  in the numerical solution of a system of ordinary differential equations}, in
  Proceedings from Symposium no. 25 held in Thessaloniki, August 17-22, 1964,
  Georgios~Ioannou Kontopoulos, ed., vol.~25 of The Theory of Orbits in the
  Solar System and in Stellar Systems, Academic Press, London, 1966, p.~281.

\end{thebibliography}
